\documentclass{article}
\usepackage[english]{babel}
\usepackage{setspace} 
\usepackage{latexsym}
\usepackage{amssymb,amsmath,amstext,amsthm,amsfonts}
\usepackage[latin1]{inputenc}
\usepackage{amsmath}
\usepackage{amsthm}
\usepackage{graphicx}
\usepackage{graphics}
\DeclareMathOperator{\Sp}{Sp}
\DeclareMathOperator{\SLc}{SL}
\DeclareMathOperator{\GL}{GL}
\DeclareMathOperator{\ot}{O}

\DeclareMathOperator{\uni}{U}

\input pictex.tex
\newcommand{\matriz}[4]{ \left(\begin{array}{cc}
#1 & #2 \\
#3 & #4\\
\end{array}\right)}

\newcommand{\C}{{\mathbb C}^\times}

\newcommand{\findem}{\qed}

\theoremstyle{definition}

\newtheorem{thm}{Theorem}[section]

\newtheorem{prop}[thm]{Proposition}

\newtheorem{lem}[thm]{Lemma}

\newtheorem{rmk}[thm]{Remark}

\newtheorem{defi}[thm]{Definition}

\newtheorem{cor}[thm]{Corollary}

\newtheorem{Proof}[thm]{Proof}

\newcommand{\nl}{\newline}
\newcommand{\SL}{\SLc^{\varepsilon}_*(2,A)}

\title{A Generalized Weil Representation for the finite split orthogonal group $\ot_q(2n,2n)$, $q$ odd greater than $3$.}
\author{Andrea Vera Gajardo.\footnote{The author was partially supported by Fondecyt Grant 1120578.}}
\date{}

\begin{document}
\maketitle

\abstract{We construct via generators and relations a generalized Weil representation for the split orthogonal group $\ot_q(2n,2n)$ over a finite field of $q$ elements.  
Besides, we give an initial decomposition of the representation found.
We also show that the constructed representation is equal to
the restriction of the Weil representation to $\ot_q(2n,2n)$ for the reductive dual pair 
$(\Sp_2(\mathbb{F}_q),\ot_q(2n,2n))$ and that the initial decomposition is the same as the decomposition with respect to the action of $\Sp_2(\mathbb{F}_q)$.\\}

\section{Introduction}

Weil representations have proven to be a powerful tool in the theory of group representations. They originate from a very general construction of A. Weil (\cite{weil}), which has as a consequence the existence of a projective representation of the group $\Sp(2n,K)$, $K$ a locally compact field. 
Weil built this representation taking advantage of the representation theory of the related Heisenberg group, as described by the Stone-von Neumann theorem in the real case (\cite{lion}).
In particular, these representations have allowed to build in a universal and uniform way all irreducible complex linear representations of  the general linear group of rank $2$ over a finite field (\cite{soto}), and later over a local field, except in residual characteristic two (\cite{ku}). 

Weil representations can be constructed in various ways. For instance, they can be constructed via Heisenberg groups, via constructions of equivariant vector bundles (\cite{gpsa}), via presentations or via dual pairs (\cite{aubert}), as we shall see later. 
The method using presentations is accomplished by having a simple presentation of the group, and 
then defining linear operators on a suitable vector space which preserve the relations among the generators of the presentation. This idea was originally suggested by Cartier, and used successfully by Soto-Andrade, for symplectic groups $\Sp(2n,\mathbb{F}_q)$ (\cite{soto}) and complex irreducible representations of $\SLc(2,\mathbb{F}_q)$.  In the case of $\Sp(2n,\mathbb{F}_q)$ he considered this group as a group "$\SLc (2)$" but with entries in the matrix ring $M_n(\mathbb{F}_q)$, and thus he obtained a suitable presentation for the group. In this way, he constructed the Weil representation for the symplectic groups $\Sp(2n,\mathbb{F}_q)$. In the case of $\Sp(4,\mathbb{F}_q)$, Soto-Andrade obtained all the irreducible representations of this group by decomposing the two Weil representations associated to the two isomorphy types of quadratic forms of rank $4$ over $\mathbb{F}_q$.

This point of view was generalized and gave rise to the groups $\SL$ for $(A,*)$ an involution ring and $\varepsilon =\pm 1 \in A$.  These groups are a generalization of special linear groups $\SLc(2, K)$, where $K$ is a field. They were defined for $\varepsilon = -1$ by Pantoja and Soto-Andrade in \cite{jofal} and generalized to $\varepsilon = 1$ in \cite{comm}. 

The groups $\SL$ include, among others, symplectic and split orthogonal groups. For instance, if $K$ is a field, $A = M_n(K)$, $\diamond$ the transposition of matrices then $\SLc_{\diamond}^{\varepsilon = -1}(2,A)$ is the symplectic group $\Sp(2n,K)$ and $\SLc_{\diamond}^{\varepsilon = 1}(2,A)$ is the split orthogonal group $\ot_K(2n,2n)$. Thus, these groups allow us to look at higher rank classical groups as rank two groups, considering them with coefficients in a new ring.

The procedure used by Soto-Andrade in \cite{soto} was approached even in the case of a non semi-simple involutive ring $(A, *)$ with non-trivial nilpotent Jacobson radical. In fact, in \cite{lucho} Guti\'errez found a Bruhat presentation and a generalized Weil representation for $\SLc_*^{\varepsilon=-1}(2,A_m)$, where $A_m=\mathbb{F}_q[x]/ \left\langle x^m \right\rangle $. 

In \cite{ongen} Guti\'errez, Pantoja and Soto-Andrade build Weil representations in a very general way via generators and relations, for the groups $\SL$ for which a ``Bruhat" presentation analogue to the classical one holds. In this way, in order to use this method it is very important to have an adequate presentation of the group. In \cite{manusc}, Pantoja generalized the classical Bruhat presentation of $\SLc(2,K)$ to $\SLc_*^{\varepsilon=-1}(2, A)$, when A is simple artinian ring with involution. 

In this work, we construct a generalized Weil representation of finite split orthogonal groups $\ot_{q}(2n,2n)$, using the method described in \cite{ongen}. As mentioned above, this group is naturally the group $\SLc_{\diamond}^{\varepsilon=1}(2,M_{2n}(\mathbb{F}_q))$. However one of the results of this paper is to realize this group as a $\SLc_{\sim}^{\varepsilon=-1}(2,M_{2n}(\mathbb{F}_q))$ group, where $\sim$ is a certain specific involution in $M_{2n}(\mathbb{F}_q)$, different from $\diamond$. This allows to use the Bruhat-like presentation that is exhibited in \cite{manusc} and facilitate significantly technical aspects of the construction. In fact, the result is more general, and provides an isomorphism between the groups $\SLc_*^+(2,M_2(A_0))$ and $\SLc_{\sim}^-(2,M_2(A_0))$, where $(A_0,*)$ is a unitary involutive ring and $\sim$ is another involution in $A_0$ obtained from $*$.

Also, we study the structure of the associated unitary group and using this group we get an initial decomposition of the representation.

In section 6 we show compatibility of the method of Guti\'errez, Pantoja and Soto-Andrade with theory of dual pairs. For this, we prove that the representation of $\ot_q(2n,2n)$ constructed using this method is equal to the restriction of the Weil representation to $\ot_q(2n,2n)$ for the dual pair $(\Sp(2,k),\ot_q(2n,2n))$. Also, we prove that the initial decomposition mentioned above is the same as decomposition with respect to the action of $\Sp(2,k)$ via the Weil representation.

\section{The groups $\SLc_{*}^{\varepsilon}(2,A)$.}

Let $(A,*)$ be a unitary ring with an involution $*$. We can extend the involution $*$ in $A$ to the ring $M_2(A)$ putting
\[
T^*=\matriz{a}{b}{c}{d} ^*=\matriz{a^*}{c^*}{b^*}{d^*}.
\]

We consider $\varepsilon = \pm 1 \in A$ and $J_{\varepsilon}= \left(\begin{array}{cc}
0 & 1 \\
\varepsilon & 0
\end{array} \right)\in M_2(A).$ Let us denote by $H_{\varepsilon}$  the associated $\varepsilon$- hermitian form defined by the matrix $J_{\varepsilon}$.

\begin{defi}
\label{defSL}
The group $\SLc_*^{\varepsilon}(2,A)$ is the set of all automorphisms
$g$ of the $A$-module $M = A \times A$ such that $H_{\varepsilon} \circ (g \times g) = H_{\varepsilon}$.
 In matrix form:

$$\SLc_*^{\varepsilon}(2,A)=\lbrace T \in M_2(A) \mid T J_{\varepsilon}T^*=J_{\varepsilon}\rbrace$$
 
\end{defi}

\begin{rmk}
In \cite{comm} it is shown that a matrix $\matriz{a}{b}{c}{d} \in M_2(A)$ is in $\SL$ if and only if the following equalities hold:
\begin{enumerate}
\item $ab^*=-\varepsilon ba^*$; 
\item $cd^*=- \varepsilon dc^*$; 
\item $a^*c=- \varepsilon c^*a$; 
\item $b^*d= - \varepsilon d^*b$; 
\item $ad^*+ \varepsilon bc^* = a^*d+ \varepsilon c^*b=1$.
\end{enumerate}

\end{rmk}

We note that if $(A,\diamond)$ is the matrix ring $M_m(\mathbb{F}_q)$ with the transpose involution $\diamond$, then 
$\SLc_{\diamond}^{-1}(2,A)$ is the symplectic group $\Sp(2m,\mathbb{F}_q)$ defined over $\mathbb{F}_q$. On the other hand 
$\SLc_{\diamond}^1(2,A)$ gives the split orthogonal group $\ot_q(m,m)$.

In what follows we put $\SLc_{*}^+(2,A)=\SLc_*^{1}(2,A)$, $\SLc_*^-(2,A)=\SLc_*^{-1}(2,A)$ and $\diamond$ always will be the transpose involution in a matrix ring.


Let $(A_0,*)$ be a unitary ring with involution $*$ and $A=M_2(A_0)$. 
Let us consider the matrix 
\mbox{ $J= \left( \begin{array}{cc}
0 & 1 \\
-1 & 0\\
\end{array}
\right) \in A^{\times}$}, which satisfies $J^{-1}=J^*=-J$. Using the matrix $J$ we can define a new involution in $A$, namely, $a^{\sim}=Ja^*J^{-1}$.

Let us consider $M_2(A)$ provided with the involutions $*$ and $\sim$ inherited from A.\\

\begin{thm}
\label{duality}
The groups $\SLc_{*}^{+}(2,A)$ and $\SLc_{\sim}^{-}(2,A)$ are isomorphic.
\end{thm}

\begin{Proof}  Let $U= \left( \begin{array}{cc}
J & 0 \\
0 & J\\
\end{array}
\right)  \in \GL_2(A) .$
A direct computation proves that $T^{\sim}U=UT^{*}$ for all $T \in M_2(A).$
It is clear that $TJ_-T^{\sim}=J_-$ if and only if $TJ_-UT^*=J_-U$. Then we must show that $J_-U$ and $J_+$ are equivalent. In fact, the (orthogonal) matrix $P=\left( \begin{array}{cc}
0 & J \\
1 & 0 \\
\end{array}
\right)\in M_2(A)$ satisfies $PJ_+P^{*}=J_-U $.
\end{Proof}

Although the split orthogonal group is naturally a ``$\SLc^+$''- group, in practical terms it is better to look at it as a ``$\SLc^-$''-group, because this fact will greatly facilitate technical aspects.\\

\begin{cor}
The split orthogonal group $\ot_q(2n,2n)$ is isomorphic to the group $\SLc_{\sim}^{-}(2,M_{2n}(\mathbb{F}_q))$, where the involution $\sim$ in $M_{2n}(\mathbb{F}_q)$is given by $a^{\sim}=J_{2n}a^{\diamond}J^{-1}_{2n}.$  ($J_{2n}= \left( \begin{array}{cc}
0 & I_n \\
-I_{n} & 0\\
\end{array}
\right) \in M_{2n}(\mathbb{F}_q)  ).$
\label{corolario}
\end{cor}

\begin{Proof} Taking $(A_0,*)$ as the involutive ring $(M_n(\mathbb{F}_q), \diamond)$ and using the theorem (\ref{duality}) we get that the groups $\SLc_{\diamond}^{+}(2,M_{2n}(\mathbb{F}_q))$ and $\SLc_{\sim}^{-}(2,M_{2n}(\mathbb{F}_q))$ are isomorphic. 

\end{Proof}

\section{A Bruhat presentation for $\SL$  }

Let $A$ be a unitary ring with involution $*$ .We will write $A^{s}_{\varepsilon,*}$ to denote the set of all $\varepsilon$- symmetric elements in A respect to the involution $*$. Namely; 
 $$A^{s}_{\varepsilon,*} = \lbrace a \in A \mid a^*=-\varepsilon a  \rbrace.$$  
In order to facilitate the notation, we put $A^{s}_{+,*}=A^{s}_{1,*}$ and $A^{s}_{-,*}=A^{s}_{-1,*}$.\\ Let us consider
\[
h_t=\left(\begin{array}{cc}
t & 0\\
0 & t^{* -1}\end{array} \right) (t \in A^{\times}),\quad  w= \left(\begin{array}{cc}
0 & 1\\
\varepsilon & 0\end{array} \right),\quad u_s=\left(\begin{array}{cc}
1 & s \\
0 & 1 \end{array} \right) (s \in A^{s}_{\varepsilon,*})
\]
\vspace{0.1cm}
\begin{defi}
\label{bruhatpres}
We will say that ${\SL}$ has a Bruhat presentation if it is
generated by the above elements with defining relations:
\begin{enumerate}
\item $h_t h_{t'}=h_{tt'}$, $u_s u_{s'}=u_{s+s'}$;
\item $w^2= h_{\varepsilon}$; 
\item$h_t u_s= u_{tst^*}h_t$;
\item $wh_t=h_{t^{*-1}}w$;
\item$wu_{t^{-1}}wu_{-\varepsilon t}wu_{t^{-1}}=h_{-\varepsilon t},\quad t \in A^{\times}\cap A^{s}_{\varepsilon,*}.$
\end{enumerate}
\end{defi}

\begin{rmk}
Observing the last relation, we note that in order to have a Bruhat presentation for ${\SL}$ is necessary that $A^{\times}\cap A^{s}_{\varepsilon,*} \neq \emptyset.$\\
\end{rmk}

In \cite{manusc} it is proved that if $A$ is a simple artinian ring with involution $*$ that either is infinite or isomorphic to the full matrix ring over $\mathbb{F}_q$ with $q>3$, then the group $\SLc_*^-(2,A)$ has a Bruhat presentation. 

Thus, the group $\SLc_{\sim}^-(2,M_{2n}(\mathbb{F}_q))$ mentioned in Corollary (\ref{corolario}) has a Bruhat presentation if $q>3$. We will use this fact to construct the desired representation.

\section{A generalized Weil representation for the split orthogonal group $\ot_q(2n,2n)$. \label{rep}}

In this section, our aim is to construct a Weil representation for the split orthogonal group seen as the group $\SLc_{\sim}^{-}(2,M_{2n}(\mathbb{F}_q))$. One way to this goal is to construct the representation using the Bruhat presentation.
For this purpose we will use the result that follows (\cite{ongen}).


Let $A$ be a ring with an involution $*$. 
Let us suppose that the ring $A$ is finite and that the group $G= \SLc_*^{\varepsilon}(2, A)$ has a Bruhat presentation.
Let $M$ be a finite right $A$-module and let us consider the following data:

\begin{enumerate}
\label{genweilrep}
\item A bi-additive function $\chi: M \times M \longrightarrow \mathbb{C}^{\times}$ and a character $\alpha \in \widehat{A}^{\times}$ such that for all $x,y \in M, t \in A^{\times}$:
\begin{enumerate}
\item $\chi(xt,y)=\alpha(tt^*)\chi(x,yt^*)$
\item $\chi(y,x)=\chi(- \varepsilon x,y)$
\item $\chi(x,y)=1$ for all $x \in M \Rightarrow y=0$
\end{enumerate}

\item A function $\gamma: A^{s}_{\varepsilon,*} \times M \longrightarrow  \mathbb{C}^{\times}$ such that for all 
$s,s' \in A^{s}_{\varepsilon,*}$, $x,z \in M$, $r\in A^{\times}$, $t \in A^{s}_{\varepsilon,*}\cap A^{\times} $ :
\begin{enumerate}
\item $\gamma(s+s',x)=\gamma(s,x)\gamma(s',x)$
\item $\gamma(s,xr)=\gamma(rsr^*,x)$
\item $\gamma(t,x+z)=\gamma(t,x)\gamma(t,z)\chi(x,zt)$
\end{enumerate}

\item $c\in \mathbb{C}^{\times}$ such that $c^2 \vert M \vert=\alpha(\varepsilon)$, 
and for all $t \in A^{s}_{\varepsilon,*} \cap A^{\times}$ the following equality holds:\\
$$\sum_{y \in M}\gamma(t,y)=\dfrac{\alpha(\varepsilon t)}{c}$$
\end{enumerate}
\begin{thm}
{\em (Guti\'errez, Pantoja and Soto-Andrade, \cite{ongen})} 
Let M be a finite right $A$-module. Denote $L^2(M)$ the vector space of all complex-valued functions on $M$, endowed with the inner product with respect to the counting measure on $M$. Set:
\begin{enumerate}
\item $\rho(h_t)(f)(x)=\overline{\alpha}(t)f(xt)$, $f\in L^2(M), t  \in A^{\times}, x \in M$;
\item $\rho(u_s)(f)(x)=\gamma(s,x)f(x)$, $f\in L^2(M), b \in A^{s}_{\varepsilon,*}, x \in M$;
\item $\rho(w)(f)(x)=c\sum_{y \in M}\chi(-\varepsilon x, y)f(y)$, $f\in L^2(M), x \in M$
\end{enumerate}

(where $\overline{\alpha}$ denotes the complex conjugate of the character $\alpha$).
These formulas define a unitary linear representation $(L^2(M), \rho)$ of $G$, called the generalized Weil representation of $G$ associated to the data $(M, \alpha, \gamma, \chi)$.

\label{gen}
\end{thm}
\vspace{0.3cm}
\begin{rmk}
Let us note that this definition contains the classic Weil representation of $\SLc_2(K)$, where $K$ is a field (see \cite{soto}, for instance).

\end{rmk}

\vspace{0.3cm}

In what follows, we will focus on finding the necessary data to construct a generalized Weil representation for $\ot_q(2n,2n)$.

From now on we put $k=\mathbb{F}_q$, $A=M_{2n}(k)$ and we consider $q$ odd greater than $3$.

We will apply Theorem (\ref{gen}) to the group $\SLc_{\sim}^{-}(2,A)\cong \ot_q(2n,2n)$.
To do this, we recall the following fact. Let $E$ be a finite dimensional vector space over a field $K$. In \cite{tignol}, the authors describe a correspondence between the linear anti-automorphisms of $End_K(E)$ and the equivalence classes of non degenerate bilinear forms on $E$ modulo multiplication by a factor in $K^{\times}$. Under this correspondence, $K$-linear involutions on $End_K(E )$ correspond to non degenerate bilinear forms which are either symmetric or skew-symmetric.
Let $B$ be a non degenerate bilinear form. The aforementioned correspondence associates $B$ with $\sigma_B$, where $\sigma_B$ is a linear anti-automorphism in $End_K(E)$ defined by the following equality;
\begin{equation}
B(f(x),y)=B(x,\sigma_B(f) y), \quad \quad  f \in End_K(E), x,y \in E.
\label{correspondencia}
\end{equation}
Now, let $V$ a vector space of $k$-dimension $2n$. We fix a basis for $V$ in order to put $M_{2n}(k)\simeq End_k(V).$ \nl
Let $<\quad,\quad>: V \times V \longrightarrow k$ be the non degenerate symmetric bilinear form given by the standard dot product. We consider the non degenerate skew-symmetric bilinear form
$[\quad,\quad]:V \times V \longrightarrow k$, given by 
$$[x,y]=<x,yJ_{2n}>.$$ 
According to the correspondence between involutions and non degenerate bilinear forms described above, the symmetric bilinear form $<\quad,\quad>$ corresponds to the transpose involution $\diamond$. Similarly, the skew-symmetric bilinear form $[\quad,\quad]$ corresponds to the new involution $\sim.$ That is, for all $x,y \in V$ and $a \in A$;

\begin{equation}
 <xa,y>=<x,ya^{\diamond}> 
\end{equation}
\begin{equation}
 [xa,y]=[x,ya^{\sim}]
 \label{invsimplec}
\end{equation}

Now, let $\psi$ be a non trivial character of $k^+$. Using the notation above let us consider:

 \begin{enumerate}
 \item 
 $M$ the right $A-$module  $V^{2}$ with the following action:
 $$(x,y)a=(xa,ya) \quad \quad a\in A ,x,y \in V.$$
 \item
 $\chi:M \times M \longrightarrow {\C},$ 
$\chi((x,y),(v,z))= \psi([x,z]-[y,v]).$
\item
$\alpha$ the trivial character of $A^{\times}.$
\item
$\gamma: A^{s}_{-,\sim}\times M \longrightarrow {{\mathbb{C}}%
^\times},$ 
$\gamma(u,(x,y))=\psi([xu,y]).$\\
\end{enumerate}

\begin{lem} 
For all $u\in A^{\times } \cap A^{s}_{-,\sim}$, the map $Q_u:V^{2}\longrightarrow k$ given by \newline
 $Q_u((x,y))=[xu,y]$ is a non degenerate split quadratic form.
 Furthermore, for $u, u^{\prime} \in A^{\times } \cap A^{s}_{-,\sim} $ the quadratic forms $Q_u$ and $Q_{u^{\prime}}$ are equivalent.

 \label{lemafc}
 \end{lem}
 
 \begin{Proof} 
Let $\lambda \in k, (x,y),(v,z)\in V^2$. 
Clearly $Q_u(\lambda (x,y))=\lambda^{2}Q_u((x,y))$. 
We will prove that
$$B((x,y),(v,z))=Q_u((x+v,y+z))-Q_u((x,y))-Q_u((v,z))$$ 
is a symmetric non degenerate bilinear form.
We have;
$$B((x,y),(v,z))= [xu+vu,y+z]-[xu,y]-[vu,z]=[xu,z]+[vu,y].$$
Now;
\begin{align*}
B((x,y)+(r,t),(v,z))=& [(x+r)u,z]+[vu,(y+t)]\\
=&[xu,z]+[ru,z]+[vu,y]+[vu,t]\\
=&B((x,y),(v,z))+B((r,t),(v,z))\\
\\
B(\lambda(x,y),(v,z))=&[\lambda xu,z]+[vu,\lambda y]\\
=&\lambda[xu,z]+\lambda [vu,y]\\
=&\lambda B((x,y),(v,z)).
\end{align*}

Then, $B$ is a symmetric bilinear form.

Let us suppose that $B((x,y),(v,z))=0$ for all $(v,z)\in V^2$. If we choose $v=0$, then $[xu,z]=0$ for all $z \in V$. Since $[\quad,\quad]$ is non degenerate and $u$ is invertible, we get $x=0$. Similarly $y=0$. Therefore, $B$ is non degenerate.

Now, if $u,u^{\prime}\in A^{\times} \cap A^{s}_{-,\sim}$ then
 $uJ_{2n}^{\diamond}$ and $u^{\prime}J_{2n}^{\diamond}$ are invertible skew symmetric matrices. In fact if $u \in A^{s}_{-,\sim} $ then $u^{\sim}=J_{2n}u^{\diamond}J_{2n}^{-1}=u$. Also $J_{2n}^{\diamond}=-J_{2n}$, so we get that $(uJ_{2n}^{\diamond})^{\diamond}=J_{2n}u^{\diamond}=uJ_{2n} =-uJ_{2n}^{\diamond}$.
 Thus $uJ_{2n}^{\diamond}$ and $u^{\prime}J_{2n}^{\diamond}$ represent a non degenerate skew symmetric bilinear form, therefore they are equivalent. So, there exists $j \in A^{\times}$ such that $uJ_{2n}^{\diamond}=ju^{\prime}J_{2n}^{\diamond}j^{\diamond}$. Thus, 
\begin{align*}
Q_{u^{\prime}}((xj,yj))&=[xju^{\prime},yj]\\
&=Q_u((x,y))
\end{align*}

If we choose $u=I_{2n}$, the quadratic form $Q_u$ is represented by the matrix $\left(\begin{array}{cc}
0 & -J_{2n}\\
J_{2n} & 0\end{array} \right)$. Thus, $Q_u$ is split.
\end{Proof}

 \begin{thm}
 
The data $(M,\alpha,\gamma,\chi)$ describe a Generalized Weil Representation for $G=\ot_q(2n,2n)$. Furthermore, this representation is independent of the choice of the character $\psi$.
\label{weilreport}
\end{thm} 
\begin{Proof} 
We will check that $\chi$ satisfies the corresponding conditions.\\ Let $(x,y), (v,z) \in M, a\in A.$
\begin{align*}
\textit{(a)} \quad \chi ((x,y)a,(v,z))=& \psi
([xa,z]-[ya,v]) \\
=&\psi([x,za^{\sim}]-[y,va^{\sim}]) \qquad ,\textit{by (\ref{invsimplec}) }\\
=& \chi ((x,y),(v,z)a^{\sim}).\\
\\
\textit{(b)}\quad \chi ((v,z),(x,y))=& \psi
([v,y]-[z,x]) \\
=&\psi([x,z]-[y,v])\\
=&\chi((x,y),(v,z)).
\end{align*}%

\textit{(c)}  Let us suppose that $\chi ((x,y),(v,z))=1$ for all $(v,z) \in M$. If $v=0$, then $\psi([x,z])=1$ for all $z\in V$. If $x \neq 0$, then $[x
,\; \cdot \; ]:V \longrightarrow k$ is a non trivial linear functional. Therefore it is surjective.
Let $\lambda \in k$ such that $\psi (\lambda)\neq 1$, and $t=t(\lambda)\in V$ such that
$\lambda=[x,t]$, then we get the following contradiction: 
\begin{equation*}
1=\psi ([x,t])=\psi (\lambda).
\end{equation*}

Therefore $x=0$, and similarly $y=0$.
\newpage
Now, we will prove that $\gamma$ satisfies the corresponding properties.
Let $u,u^{\prime }\in A^{s}_{-,\sim}$, $a\in A^{\times}$, $(x,y),(v,z) \in M;$

\begin{align*}
\textit{(a)}\quad  \gamma(u+u^{\prime},(x,y)) 
=&\psi([xu+xu^{\prime},y]) \\
=&\psi([xu,y])\psi([xu^{\prime},y])\\
=&\gamma(u,(x,y)) \gamma(u^{\prime },(x,y)).\\
\\
\textit{(b)} \quad \gamma(u,(x,y)a)=& \psi([xau,ya]) \\
=&\psi([xaua^{\sim},y]) \\
=&\gamma(aua^{\sim},(x,y)).\\
\\
\textit{(c)} \quad \gamma(u,(x,y)+(v,z))=&\psi([(x+v)u,y+z]) \\
=&\psi([xu,y])\psi([xu,z])\psi([vu,y])\psi([vu,z]) \\
=&\gamma(u,(x,y))\gamma(u,(v,z))\psi([x,zu]-[y,vu]) \\
=& \gamma(u,(x,y))\gamma(u,(v,z))\chi((x,y),(v,z)u).
\end{align*}

Now, we must choose $c \in {\C}$ satisfying $c^{2}|M|=1$ and show that for $u\in A^{\times }\cap A^{s}_{-,\sim}$ the following equality holds:

\begin{equation*}
\sum_{(x,y)\in M}\gamma (u,(x,y))=\sum_{x,y\in V}\psi ([xu,y])=\dfrac{1}{c}.
\end{equation*}

 According to the lemma \ref{lemafc}, we know that $\sum_{(x,y) \in M}\gamma(u,(x,y))$ is a Gauss sum associated to a split quadratic form in a vector space of even dimension $4n$. This sum is calculated, for instance, in \cite{soto}.
 In fact, $$\sum_{(x,y) \in M}\gamma(u,(x,y)) = q^{2n},$$

 We choose $c=\dfrac{1}{q^{2n}}$. Thus,
$$ \sum_{(x,y)\in M}\gamma(u,(x,y))=\dfrac{1}{c}.$$
 
 Now, let $\psi_1$ and $\psi_2$ be two non trivial characters of $k^+$. Let us prove that the corresponding representations are isomorphic.

 Let $\lambda \in k^{\times}$ such that $\psi_2(r)=\psi_1(\lambda r)$ for all $r \in k$. 
 Let $(L^2(M),\rho_1)$ and $(L^2(M), \rho_2)$ the Weil representations obtained from $\psi_1$ and $\psi_2$ respectively. Then, the linear automorphism \newline
 $\Psi:L^2(M) \longrightarrow L^2(M)$ given by $(\Psi f)(x,y)=f(x,\lambda y)$ is a isomorphism between the representations $(L^2(M),\rho_1)$ and $(L^2(M), \rho_2)$.
 
 \end{Proof}
 
\section{An initial decomposition.}
\label{decomp}

\begin{defi}
\label{unit}
The group $\uni(\gamma, \chi)$ is the group of all $A$-linear automorphisms $\beta$ of $M$ such that:
\begin{enumerate}
\item $\gamma(u,\beta(x,y))=\gamma(u,(x,y))$ for all $u \in A^{s}_{\varepsilon,*},(x,y)\in M.$
\item $\chi(\beta(x,y),\beta(v,z))=\chi((x,y),(v,z))$ for all $(x,y),(v,z) \in M.$
\end{enumerate} 
In what follows we will denote $\uni(\gamma, \chi)$ simply by $\uni$.
\end{defi}
 
 \newpage
Following the idea of \cite{ongen}, if we know the structure of the group $\uni$ and the set of its irreducible representations, we can find an \textit{initial} decomposition of the Weil Representation in the sense that we do not know if the components obtained are irreducible. In what follows, we make this decomposition explicit.

For $\beta \in \uni$ and $x \in M$ we put $ \beta.x=\beta(x)$.
 The group $\uni$ acts naturally on $L^2(M)$. That is to say  the action is given by:

\begin{align*}
&\sigma: \uni \longrightarrow Aut_{\mathbb{C}}(L^2(M)),\\
&\sigma_{\beta}(f)(x)=f(\beta^{-1}.x)
\end{align*}

In \cite{ongen} it is shown that the natural action of $\uni$ on $L^2(M)$ commutes with the action of the Weil Representation.

\vspace{0.2cm}

Let $\widehat{\uni}$ be the set of the irreducible representations of $\uni$. 
We consider the isotypic decomposition of $L^2(M)$ with respect to $\uni$:
\[
L^2(M) \cong \bigoplus_{(V_{\pi},\pi) \in \widehat{\uni}} n_{\pi}V_{\pi}.
\]
Since $n_{\pi}=dim_{\mathbb{C}}(Hom_{\uni}(V_{\pi},L^2(M))) =dim_{\mathbb{C}}(Hom_{\uni}(L^2(M),V_{\pi}))$, we can write this decomposition in the following  way:

\[
L^2(M) \cong \bigoplus_{(V_{\pi},\pi) \in \widehat{\uni}} (Hom_{\uni}(L^2(M),V_{\pi})\otimes_{\mathbb{C}} V_{\pi}.)
\]
If we put $m_{\pi}=dim_{\mathbb{C}}(V_{\pi})$, we get;
\[
L^2(M) \cong \bigoplus_{(V_{\pi},\pi) \in \widehat{\uni}} m_{\pi}Hom_{\uni}(L^2(M),V_{\pi}.)
\]
If $(V_{\pi},\pi)\in \widehat{\uni}$ and $\beta \in \uni$, we denote by $\pi_{\beta }$ the map $\pi(\beta):V_{\pi} \longrightarrow V_{\pi}$.
The space $Hom_{\uni}(L^2(M),V_{\pi})$ is formed by linear functions $\Theta:L^2(M) \longrightarrow V_{\pi}$ such that for any $\beta \in \uni$  
\begin{equation}
\Theta \circ \sigma_{\beta} =\pi_{\beta} \circ \Theta .
\label{entre}
\end{equation}

Let us consider the Delta functions $\lbrace e_x \mid x \in M \rbrace$ and the map  $\theta:M \longrightarrow V_{\pi}$  such that $\theta(x)=\Theta(e_x)$ for all $x \in M$. Since $\sigma_{\beta}(e_x)=e_{\beta.x}$, condition (\ref{entre}) becomes:

\begin{equation}
\theta(\beta.x)=\pi_{\beta} \circ \theta(x).
\label{casi}
\end{equation}
Conversely, let $\theta:M \longrightarrow V_{\pi}$ satisfying (\ref{casi}). We extend linerarly and we get a map $\Theta:L^2(M) \longrightarrow V_{\pi}$ such that (\ref{entre}) holds.

Thus, we can see the space $Hom_{\uni}(L^2(M),V_{\pi})$ as the function space formed by maps
$\theta:M \longrightarrow V_{\pi}$ such that $\theta(\beta.x)=\pi_{\beta}\circ \theta(x)$ for all $\beta \in \uni, x \in M$. The group $G=\SL$ acts on this space via the Weil representation, using the same formulas as defined in Theorem (\ref{gen}). Similarly, it is possible to define the natural action of the group $\uni$ in this space, because- like $L^2(M)$- it is formed for functions with domain $M$.
 
\newpage 
Let $\rho$ denote the Weil action of $G$ on $L^2(M)$ and $\widehat{\rho}$ the Weil action of $G$ on $\bigoplus_{(V_{\pi},\pi) \in \widehat{\uni}} m_{\pi}Hom_{\uni}(L^2(M),V_{\pi})$.
Because of how we define the Weil representation, there exist scalars $K_g(x,y)\in \mathbb{C}$ depending only on $g \in G$ and $x,y \in M$ such that for all $f \in L^2(M)$, $\Lambda \in \bigoplus_{(V_{\pi},\pi) \in \widehat{\uni}} m_{\pi}Hom_{\uni}(L^2(M),V_{\pi})$ the following statements holds:

\begin{align}
\rho_g(f)&= \sum_{y \in M}K_g(\cdot,y)f(y);\\
\widehat{\rho}_g(\Lambda)&=\sum_{y \in M}K_g(\cdot,y)\Lambda(y). \label{a}
\end{align}
In this way, we get:\\

\begin{lem}
 $(L^2(M),\rho)$ and $(\bigoplus_{(V_{\pi},\pi) \in \widehat{\uni}} m_{\pi}Hom_{\uni}(L^2(M),V_{\pi}),\widehat{\rho})$ are isomorphic representations of $G$.
\end{lem}

\begin{Proof}
The linear isomorphism between $L^2(M)$ and $\bigoplus_{(V_{\pi},\pi) \in \widehat{\uni}} m_{\pi}Hom_{\uni}(L^2(M),V_{\pi})$
is an isomorphism between representations.
\end{Proof}

Finally, we have:\\

\begin{prop}
{\em The space $Hom_{\uni}(L^2(M),V_{\pi})$ is invariant under the Weil action of $G$.}
\end{prop}

\begin{Proof}
 Let $g \in G$, $\theta \in Hom_{\uni}(L^2(M),V_{\pi})$, $\beta \in \uni$, $x\in M$:
\begin{align*}
(\widehat{\rho}_g \theta) (\beta.x)&= \sigma_{\beta^{-1}}(\widehat{\rho}_g \theta)(x),\quad \textit{(definition of } \sigma_{\beta})\\
&=\widehat{\rho}_g(\sigma_{\beta^{-1}}\theta)(x)\\
&=\widehat{\rho}_g(\pi_{\beta} \circ \theta)(x), \quad (\sigma_{\beta^{-1}}(\theta)=\pi_{\beta} \circ \theta)\\
&= \pi_{\beta}(\widehat{\rho}_g \theta(x).)
\end{align*}
 The last equality holds because (\ref{a}).
\end{Proof}
\vspace{0.5cm}

Now, having made the decomposition above explicit, our purpose is to obtain an initial decomposition for our particular case $G=\ot_q(2n,2n) \cong \SLc_{\sim}^{-}(2,M_{2n}(k))$. For this it is enough to know the structure of the group $\uni$ and the set of irreducible representations.\\

\begin{rmk}
We note that since in our case $A=M_{2n}(k)$ and $A^{s}_{-,\sim}\cap A^{\times} \neq \emptyset $, the first condition in definition (\ref{unit}) implies the second one (see \cite{ongen}).\\
\end{rmk}

\begin{thm}
Let $\gamma$ and $\chi$ be the functions defined above. Then, $$\uni(\gamma,\chi)\cong \SLc_2(k)$$.
\end{thm}

\noindent{\em Proof.} Let $\beta \in \uni(\gamma, \chi)$. In particular $\beta$ is  $k-$ linear, therefore we can suppose that $\beta \in M_{4n}(k)$. We can write the action of $A$ on $M$ in matrix language as follows:
 
 \begin{equation*}
  (x \quad y) \matriz{a}{0}{0}{a}=(xa \quad ya)\qquad x,y \in V, a \in A.
 \end{equation*}

Since $\beta$ is $A-$linear we have that $\beta(x,y)a=\beta(xa,ya)$. In matrix language;
\begin{equation}
\beta \matriz{a}{0}{0}{a}=\matriz{a}{0}{0}{a} \beta.
\label{eq}
\end{equation}

Let $\beta_1, \beta_2, \beta_3, \beta_4 \in A$ such that $\beta = \matriz{\beta_1}{\beta_2}{\beta_3}{\beta_4}$. Then, using (\ref{eq}) we get that each of these blocks must be scalar. Thus, there are $b_1,b_2,b_3,b_4 \in k$ such that  $\beta=\matriz{b_1 I_{2n}}{b_2 I_{2n}}{b_3 I_{2n}}{b_4 I_{2n}}$ and hence
\begin{equation*}
\gamma (u,\beta(x,y))=\psi([(b_1x+b_3y)u,b_2x+b_4y]).
\end{equation*}

Let us note that the bilinear form  $(x,y)\mapsto[xu, y]$ is skew symmetric for all $u  \in A^{s}_{-,\sim}$, hence $[xu,x]=0$ for all $x \in V, u \in A^{sym}$. Thus for all $x,y \in V, u \in A^{s}_{-,\sim}$
\begin{align*}
\gamma (u,\beta(x,y))=&\psi((b_1b_4-b_2b_3)[xu,y])\\=&\psi([xu,y])\\=&\gamma (u,(x,y)).
\end{align*}

Consequently 
$\psi((b_1b_4-b_2b_3-1)[xu,y])=1$ for all $x, y \in V, u \in A^{s}_{-,\sim}.$

From this last equality it follows that $b_1b_4-b_2b_3=1$. In fact, let $u \in A^{s}_{-,\sim} \cap A^{\times}$, $x \neq 0$ and let us suppose $b_1b_4-b_2b_3-1 \neq 0$.

The map $F_{x,u}:V \longrightarrow k$ given by $F_{x,u}(z)=[(b_1b_4-b_2b_3-1)xu,z]$ is a non trivial linear functional and therefore is surjective. 
Let $\lambda \in k$ such that $\psi (\lambda) \neq 1$ and $z=z(\lambda) \in V$ such that \newline $\lambda=[(b_1b_4-b_2b_3-1)xu,z]$. Then 
$\psi(\lambda)= \psi([(b_1b_4-b_2b_3-1)xu,z])=1. $ This contradicts our assumption and therefore our result follows.
$\findem$
\vspace{0.6cm}

Thus, for our case, we get an initial decomposition of the Weil Representation $(L^2(M),\rho)$. We expect to address the question about irreducibility elsewhere.

\section{Dual Pairs.}

In this section we will prove that the representation $(L^2(M),\rho)$ of $\ot_q(2n,2n)$ constructed in section \ref{rep} is equal to the restriction of the Weil representation to $\ot_q(2n,2n)$ for the dual pair $(\Sp(2,k),\ot_q(2n,2n))$. Also, we will prove that the initial decomposition described above is the same as decomposition with respect to the action of $\Sp(2,k)$ via the Weil representation.

Let $J_{2n}= \left(\begin{array}{cc}
0 & I_n \\
-I_n & 0
\end{array} \right)\in M_{2n}(k)$ and $F= \left(\begin{array}{cc}
0 & J_{2n} \\
-J_{2n} & 0
\end{array} \right)\in M_{4n}(k).$
The matrix $F$ defines the following non-degenerate split symmetric bilinear form in $V_1=k^{4n}$ 
$$(u,v)=v^tFu \qquad u,v \in V_1$$
The group $G$ of isometries of this form is isomorphic to the split orthogonal group $\ot_q(2n,2n)$.
As before, set
$$a^{\sim}=J_{2n}a^tJ_{2n}\qquad a \in M_{2n}(k).$$

A direct calculation shows that the following matrices belong to the group $G$:
\[
h_a=\left(\begin{array}{cc}
a & 0\\
0 & (a^{\sim})^{ -1}\end{array} \right) ,\quad  w= \left(\begin{array}{cc}
0 & I_{2n}\\
-I_{2n} & 0\end{array} \right),\quad u_s=\left(\begin{array}{cc}
I_{2n} & s \\
0 & I_{2n} \end{array} \right) 
\]
$(a \in M_{2n}(k)^{\times}, s =s^{\sim} \in M_{2n}(k))$.

Therefore $G=\SLc_{\sim}^{-}(2,M_{2n}(k)$.

\newpage 
Let $V_2= k^2$ and $W=Hom(V_1,V_2)$. The following formula defines a non-degenerate symplectic form on $W$.
$$\ll w_1,w_2 \gg =tr(w_1Fw_2^tJ_2)\quad \quad (w_1,w_2 \in W )$$
The group $G$ acts on $W$ by 
$$g(w_1)=w_1g^{-1} \quad \quad (g \in G, w_1 \in W).$$
This action preserves the symplectic form $\ll\cdot,\cdot \gg$. In fact, since $g \in G$, 
$$\ll w_1g^{-1}, w_2g^{-1} \gg =tr(w_1g^{-1}F(g^{-1})^tw_2^tJ_2)=tr(w_1Fw_2^tJ_2)= \ll w_1, w_2 \gg.$$

Let 
$$X= \lbrace(x,0) \mid x \in M_{2,2n}(k) \rbrace ;\quad Y= \lbrace (0,y)\mid y \in M_{2,2n}(k) \rbrace.$$
Then $W=X\oplus Y$ is a complete polarization. We will consider the Schr\"odinger model of the Weil representation of $\Sp(W)$ attached to the above complete polarization realized on $L^2(X)$ as in \cite{aubert}. Let $(L^2(X),\omega)$ such representation.

We identify $X$ with $M_{2,2n}(k)$ in the canonical way
 $$X \ni (x,0) \leftrightsquigarrow x \in M_{2,2n}(k).$$
 \begin{rmk}
 Let us note that the module $M$ in section \ref{rep} is canonically isomorphic to $X$. Consequently the spaces $L^2(M)$ and $L^2(X)$ are also isomorphic.
 \end{rmk}
 
Let $\psi$ a non-trivial character of the additive group $k^+$.
For all $x \in X$, $y \in Y$ it is clear that $h_a(x)=xh_a^{-1} \in X$ and $h_a(y)=yh_a^{-1} \in Y$, then the matrix $h_a$ preserves $X$ and $Y$. Also, $det(h_a \vert_X)\in k^{\times 2}$. Thus, proposition 34 in \cite{aubert} shows that:
$$\omega (h_a)f(x)=f(xa)  \quad \quad (f \in L^2(X))$$
Thus, $\omega(h_a)=\rho(h_a).$

Now, let us see the action of $\omega$ on $u_s$. The matrix $u_s$ acts trivially on $Y$ and on $W/Y$. Therefore, proposition 35 in \cite{aubert} shows that:
$$\omega(u_s)f(x)=\psi(\ll xc(-u_s),x \gg)f(x),$$
where $c(-u_s)= \left(\begin{array}{cc}
0 & -s/2\\
0 & 0\end{array} \right) \in M_{4n}(k)$ is the Cayley transform for $-u_s$.

Let $x= \left(\begin{array}{cc}
x_1 & 0\\
x_2 & 0\end{array} \right)\in X $, $x_1,x_2 \in k^{2n}$. Then,
 $$\ll xc(-u_s),x \gg=x_2sJ_{2n}x_1^t.$$
 In order to prove that $\omega(u_s)=\rho(u_s)$ we have to check that 
 $$[x_1s,x_2]=\ll xc(-u_s),x \gg,$$
 where $[\quad,\quad]$ is the symplectic form defined in section $4$. In fact,
 $$[x_1s,x_2]=[x_1,x_2s]=-[x_2s,x_1]=x_2sJ_{2n}x_1^t$$
\newpage 
 It is clear that the matrix $w$ maps $X$ bijectively onto $Y$ and $Y$ onto $X$, and $w^2=-1$. Then, using proposition 36 of \cite{aubert} we get:
 $$\omega(w)f(x)=\dfrac{1}{\sqrt{\vert X \vert}}\sum_{x'\in X}\psi (\ll w(x),x' \gg)f(x')$$
 Thus, in order to prove that $\omega(w)=\rho(w)$ we have to check that 
 $$\chi(x,x')=\psi (\ll xw^{-1},x' \gg).$$
 Let $x= \left(\begin{array}{cc}
x_1 & 0\\
x_2 & 0\end{array} \right) $,$x'= \left(\begin{array}{cc}
x'_1 & 0\\
x'_2 & 0\end{array} \right) $, $x_1,x_2, x'_1, x'_2 \in k^{2n}$. So, 
$$\psi (\ll xw^{-1},x' \gg)=\psi(x_2J_{2n}(x_1')^t-x_1J_{2n}(x_2')^t).$$
On the other hand, 
\begin{align*}
\chi(x,x')&=\psi([x_1,x'_2]+[x'_1,x_2])\\
&=\psi (\left( x_1 \quad x_2 \right) \left(\begin{array}{cc}
0 & -J_{2n}\\
J_{2n} & 0\end{array} \right) \left(\begin{array}{c}
(x'_1)^t \\
(x'_2)^t \end{array} \right))\\
&=\psi(x_2J_{2n}(x_1')^t-x_1J_{2n}(x_2')^t)).
\end{align*}
Thus, we have showed that the representation constructed in section \ref{rep} is equal to the restriction of the Weil representation to $O_q(2n, 2n)$ for the dual pair $(\Sp(2,k),\ot_q(2n,2n))$.

Furthermore, since an element $g \in \Sp(2,k)=\SLc_2(k)$ preserves $X$ and $Y$ and $det(g\vert_X)\in k^{\times 2}$, using proposition 34 in \cite{aubert} we get that the group $\SLc_2(k)$ acts on $L^2(X)$ as follows:
$$\omega(g)f(x)=f(g^{-1}x) \quad g \in \SLc_2(k), f \in L^2(X), x \in X.$$

Therefore, the initial decomposition in section \ref{decomp} is the same as the decomposition with respect to the action of $\SLc_2(k)$ via the Weil representation.

\vspace{1cm}
\begin{flushleft}
\textbf{Acknowledgments:} The author would like to thank Jos\'e Pantoja, Jorge Soto, Luis Guti\'errez and Anne-Marie Aubert for their never ending support and willingness. 
\end{flushleft}

\vspace{2cm}
\begin{flushleft}
Andrea Vera Gajardo. \newline
Universidad de Santiago de Chile. \newline
Avenida Libertador Bernardo O'Higgins 3363, Estaci\'on Central, Santiago, Chile. \newline
email: andreaveragajardo@gmail.com
\end{flushleft}

\end{document}